\documentclass[runningheads]{llncs}
\usepackage{amsmath} % align
\usepackage{graphicx}
\usepackage[svgnames]{xcolor} % defines e.g. "olive" for maple.sty
\usepackage{hyperref}
\usepackage{breqn} % defines dmath environment for maple.sty
\usepackage{maple}
\usepackage{mathtools} % defines e.g. coloneqq for maple.sty
\usepackage{listings}
\definecolor{mygreen}{rgb}{0,0.6,0}
\definecolor{mygray}{rgb}{0.5,0.5,0.5}
\definecolor{mymauve}{rgb}{0.58,0,0.82}
\definecolor{altblue}{rgb}{0.0,0.6,1.0}
\definecolor{lstbg}{cmyk}{0.05, 0.01, 0, 0}
\definecolor{morebluish}{cmyk}{0.06,0.04,0,0}
\input{listings-maple-definition.sty}
\newcommand{\mat}[1]{\ensuremath{\mathbf{#1}}}
\begin{document}
\title{ Teaching Linear Algebra \\ in a mechanized mathematical environment\thanks{Supported in part by NSERC and by the MICINN.}}
\titlerunning{Teaching Linear Algebra}
% If the paper title is too long for the running head, you can set
% an abbreviated paper title here
%
\author{Robert M.~Corless\inst{1}\orcidID{0000-0003-0515-1572} \and
David J.~Jeffrey\inst{1}\orcidID{0000-0002-2161-6803} \and
Azar Shakoori\inst{2}\orcidID{0009-0003-6181-4665}}
\authorrunning{R.M.~Corless et al.}
% First names are abbreviated in the running head.
% If there are more than two authors, 'et al.' is used.
%
\institute{University of Western Ontario, London, Ontario, Canada
\email{\{rcorless,djeffrey\}@uwo.ca}\\
\and
Ontario Tech University, Oshawa, Ontario, Canada\\
\email{Azar.Shakoori@ontariotechu.ca}}
\maketitle              % typeset the header of the contribution
\begin{abstract}
This paper outlines our ideas on how to teach linear algebra in a mechanized mathematical environment, and discusses some of our reasons for thinking that this is a better way to teach linear algebra than the ``old fashioned way''.  We discuss some technological tools such as Maple, Matlab, Python, and Jupyter Notebooks, and some choices of topics that are especially suited to teaching with these tools.  The discussion is informed by our experience over the past thirty or more years teaching at various levels, especially at the University of Western Ontario.
\end{abstract}
\section{Overview}
\begin{quote}
    ``Linear algebra is the first course where the student encounters algebra, analysis, and geometry all together at once.''
\end{quote}
\begin{flushright}
    ---William (Velvel) Kahan, \\
    to RMC at the 4th SIAM Linear Algebra Conference in Minneapolis 1991
\end{flushright}
This paper describes the current state of our ongoing practice of teaching linear algebra in mechanized environments. We report our thoughts, arrived at after several decades of history in differing technological and administrative support structures.  Some of our teaching philosophy is laid out in~\cite{Betteridge2022} and the references therein (especially for \textsl{active} teaching), but to keep this paper self-contained we will give a precis of our approach in section~\ref{sec:Approach}.

We believe that this paper will be of interest for this conference both for its use of various computational environments (Jupyter notebooks, Maple, Matlab, and historically the HP48 series of calculators) and for its recommendations of what is needed for future environments for mechanized mathematics.

Linear algebra as a mathematical subject is second only to Calculus in terms of overall teaching effort at secondary institutions, accounting for many millions of dollars spent every year. There are those who believe that we should devote even more money and effort to it, because linear algebra is foundational for so many applications: optimization (linear programming), scientific computing, and analysis of data, for examples.

We take as fundamental that the vast majority of people taking these enormous numbers of courses are \textsl{not} going to choose careers as pure mathematicians.  Rather, they are going to become engineers, biologists, chemists, physicists, economists, computer scientists, or something else\footnote{The diversity of where our students go afterwards makes it tricky to choose motivating applications.  Network flow problems will appeal to a subset of people; electrical circuits might appeal to another subset.  Markov chains are fun for some.  Very few applications are interesting to \textsl{everybody}.}.  They will likely need probability, and methods to solve linear equations, and the understanding of what an eigenvalue is (and perhaps what a singular value is).  By and large they will not need to reason their way out of tricky artificial problems.  They will need graph theory, and how to solve algebraic equations.  They will need to learn how to use computers to help with the drudgery of the computations involved, so that they can be free to think about what the answers mean, instead of how they are arrived at.  They will need to learn when they can rely on computers to help, and when they should be suspicious.

% The number, diversity, and quality of existing linear algebra books is substantial.  Almost every reader of this paper will have two or three ``favourites''.  If not, one can find a list of ten pretty good choices at \url{https://linearalgebras.com/ten-best-linear-algebra-books.html}.  Not all of those books are suitable for the large numbers of people described in the paragraph above, however.  Some are more appropriate for future mathematicians.

% Not listed at that URL, but one of our favourites even though it is not a textbook, is the CRC Handbook of Linear Algebra~\cite{hogbenhandbook} (perhaps we are biased, because two of us wrote the chapter therein on Linear Algebra in Maple).  There are at least a dozen absolutely excellent \textsl{numerical} linear algebra textbooks, as well: the classic~\cite{golub2013matrix} has been cited an astonishing number of times, and deservedly so.  We also recommend~\cite{watkins2004fundamentals} and~\cite{meyer2000matrix} as being exceptionally clear and comprehensive.  Perhaps~\cite{trefethen2022numerical} is the most accessible of the ones we list here.

Our favourite \textsl{introductory} textbook---out of the myriad possible choices---arose from an NSF-funded educational project, namely~\cite{carlson1993linear}.  The book is~\cite{lay2016linear}.  Yet this choice is not uncontroversial, and the book is not an especially good match for a mechanized environment.  We see a need for a specialized textbook to support active learning of linear algebra in a mechanized environment.
%We will try to support that assertion in what follows.
\subsection{Active learning in a mechanized environment\label{sec:Approach}}
Within the mathematics mechanization community, it is uncontroversial to assert that the tools available and being developed will make the learning and practice of mathematics better.  In theory, this is obvious.  In practice, there are devils in the details.  For one thing, students (and researchers in industrial environments) must be trained in the use of the new tools, and the time spent learning these tools cannot also be spent on learning the mathematical topics.  For this reason, we advocate at least some ``re-use'' of tools, namely that teaching of mechanized mathematics should use tools that will also be used for something else in the student's or researcher's career.

Nowadays this largely means Jupyter notebooks and Python, which are both very popular in data science and neuroscience.  In a few years this might mean a replacement for Jupyter together with Julia (perhaps).  The one thing that we can say about the software environment for mathematics is that it is changing as rapidly now as it ever has been.

However, it will not be surprising to the attendees of this conference that there are lessons to be learned from attempts to use mechanized mathematics in teaching in the past.  Indeed the ``deep structure'' of Python is not so different from that of Maple, and many aspects of programming in the one language transfer readily to the other (for instance, dictionaries in Python are analogous to tables in Maple).  More to the point, learning to program in \textsl{any} language exercises some of the same mental muscles that writing a mathematical proof does.  The analogy between recursion and mathematical induction is very close, indeed.  So, at least some of the material that has been developed with older technology can be given some syntactic re-sugaring and used in much the same way.  We will give examples.

The most important use of technology, however, is to increase the activity level of the student.  One needs to engage the student's attention, and get them to do more than just passively read a text, attend a lecture, watch a video, or regurgitate on an exam.  In some ways, fashion helps with this.  The students are more likely to want to learn Python than (say) C.

\subsection{How to teach with technology}
There are many papers, and indeed books, written on how to teach with technology.  We mention the influential paper~\cite{buchberger1990should}, which introduced the ``White Box'' / ``Black Box'' model, which we have used with some success.  The idea there is that when teaching a particular technique (for instance, what a determinant is) the student is not allowed to use the \lstinline{Determinant} command; but after they have understood that topic, whenever they are \textsl{using} determinants in a future topic (say, Cramer's Rule) they are allowed to use it.  The psychological and pedagogical point is that people need a certain amount of human action with a concept before it is internalized.  We tend to say that at that point, the concept has become an \textsl{answer} to the student instead of a \textsl{question}.  At that point, the students can use the technology with assurance, and the feeling that they know what is going on.

This rule can be used in other ways, and even backwards: use a tool as a mysterious Black Box for a while, probing its output by giving it various inputs until some sense of what is going on arises.  We have used this reverse strategy with some success, as well, most commonly with the Singular Value Decomposition (SVD).  See~\cite{Betteridge2022} for more strategies for teaching with technology that have been tested in practice.

\subsection{What to teach, when technology is involved}
A much more interesting question arises when one considers that the curriculum must be continually curated as new tools come available.  New topics may be added (for instance, the SVD), and old topics dropped (for instance, condensation, or perhaps Gauss--Seidel iteration).  Indeed a certain amount of room must be made in the course for instruction in the responsible use of the new tools.  This is by no means easy, and the students will resist such instruction if they are not also assessed on the use of the tools.  The fact that they will be expected to use these tools later in life as a matter of course is sometimes not enough to encourage the students to learn them now.  However, society appears to expect that we as instructors will be teaching the students the best way to actually use the material we teach, and (as a matter of course) this means that we must be teaching the students to use the tools of modern mechanized mathematics.  Those of us who are actually in the classroom know that sometimes compromises are necessary.

\subsection{Outline of the paper} In section~\ref{sec:tools}, we discuss some of the tools that are available. In section~\ref{sec:topics} we mention a few necessary topics that work well with these tools (we do not give a full syllabus, because of space limitations).  In section~\ref{sec:assessment} we discuss methods of assessment. In section~\ref{sec:reaction} we discuss some reactions from colleagues and students to these changes from a traditional syllabus, and then conclude.

\section{Tools\label{sec:tools}}
The members of this community will have their own preferred computational tools, which may not be the same as ours.  We will not fully justify our choices here, but instead sketch only some of the reasons for our choices.
\subsection{Proprietary Tools}
We do use some proprietary tools, namely Maple and Matlab.  Our Universities have site licences for these, and we have a significant body of experience with using these tools both for research and for teaching.  Many engineering students will graduate into work environments that have Matlab, and by the usual feedback mechanism from other students and other professors, most engineering students are well-motivated to learn Matlab.  Matlab has some especially nice tools for sparse matrices, and its live scripts are quite usable.

Maple is less well-used in industry, but in some countries it does have a presence; nonetheless it is a harder ``sell'' to students, and if the course does not explicitly give marks for knowing how to use Maple, students are sometimes reluctant to spend time learning it.  But it is powerful enough that students do appreciate it, once they have made the effort.

There are other proprietary products which also could be used.  Maple Learn is a new one, for instance; but we do not yet have experience with it.

Other places will use Mathematica instead of Maple, but the concerns and affordances are similar.
\subsection{Free software}
Within the free software ecosystem, Python and Jupyter stand out as tools of choice for a lot of scientists and engineers.  For linear algebra, Matlab and Maple are both superior in terms of capability and in terms of ease of use (in our opinion), especially for sparse matrices, but there is no doubt whatever that Python and Jupyter are more popular.

% Part of that popularity may be due to Markdown, which gives a very lightweight method for including math equations into the text of a Jupyter notebook.  Markdown is noticeably easier to learn than is \LaTeX, although Markdown is not as powerful.  Using Markdown is also easier (in our opinion) than using a palette-based approach, which we usually find painfully slow to write with.  Markdown also produces quite legible output, owing in part to the hard work of people in the mathematics mechanization community.

Python is remarkable for its support for long integer arithmetic (although its quiet casting of types behind the scenes can cause problems, especially when things unexpectedly contain 32 bit unsigned integers instead of the expected long integers).  Learning to program in Python is perhaps easier in the beginning than is learning any other language (we are aware that opinions differ in this regard, but surely the statement ``the easy parts of Python are easy to learn'' would be uncontroversial).

Julia is newer, more exciting, and extremely impressive for its speed as well as its ease of use.  We anticipate that use of Julia will eclipse that of Python.
% In particular, its tools for numerical solution of differential equations are extremely impressive.  See~\cite{rackauckas2017differentialequations}.

\subsection{Visualization}
Linear algebra might not seem to need visualization tools as much as Calculus does, but there are several instances where we have found dynamic visualizations to be extremely helpful.  One is exemplified by the old Matlab command \lstinline{eigshow} (which, curiously, has been deprecated and moved into a relatively obscure location inside the Matlab environment) which is extremely effective in giving students ``aha!'' moments about both eigenvalues and singular values.  One of the keys to that tool's effectiveness is (was) the kinesthetic use of the mouse, by the student, to move the input vectors around.  The immediate visual feedback of where the output eigenvectors (and singular vectors) move to in response is, in our experience, \textsl{much} more effective than simple animations (or static pictures).

More simply, getting the students to plot eigenvalue distributions, or to plot eigenvector components, is valuable as an action.

\begin{figure}
    \centering
    \includegraphics[width=2cm]{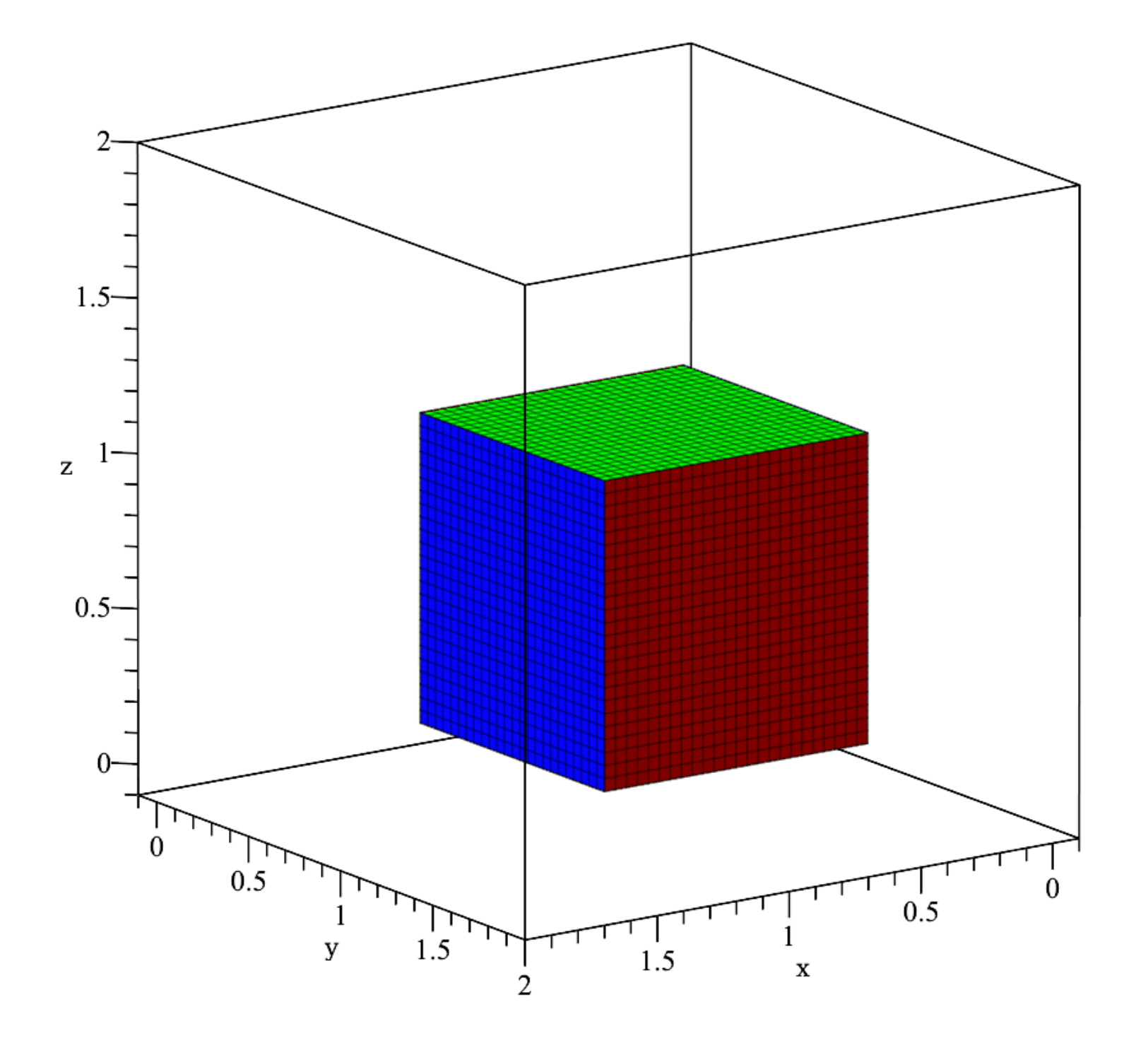}\includegraphics[width=3cm]{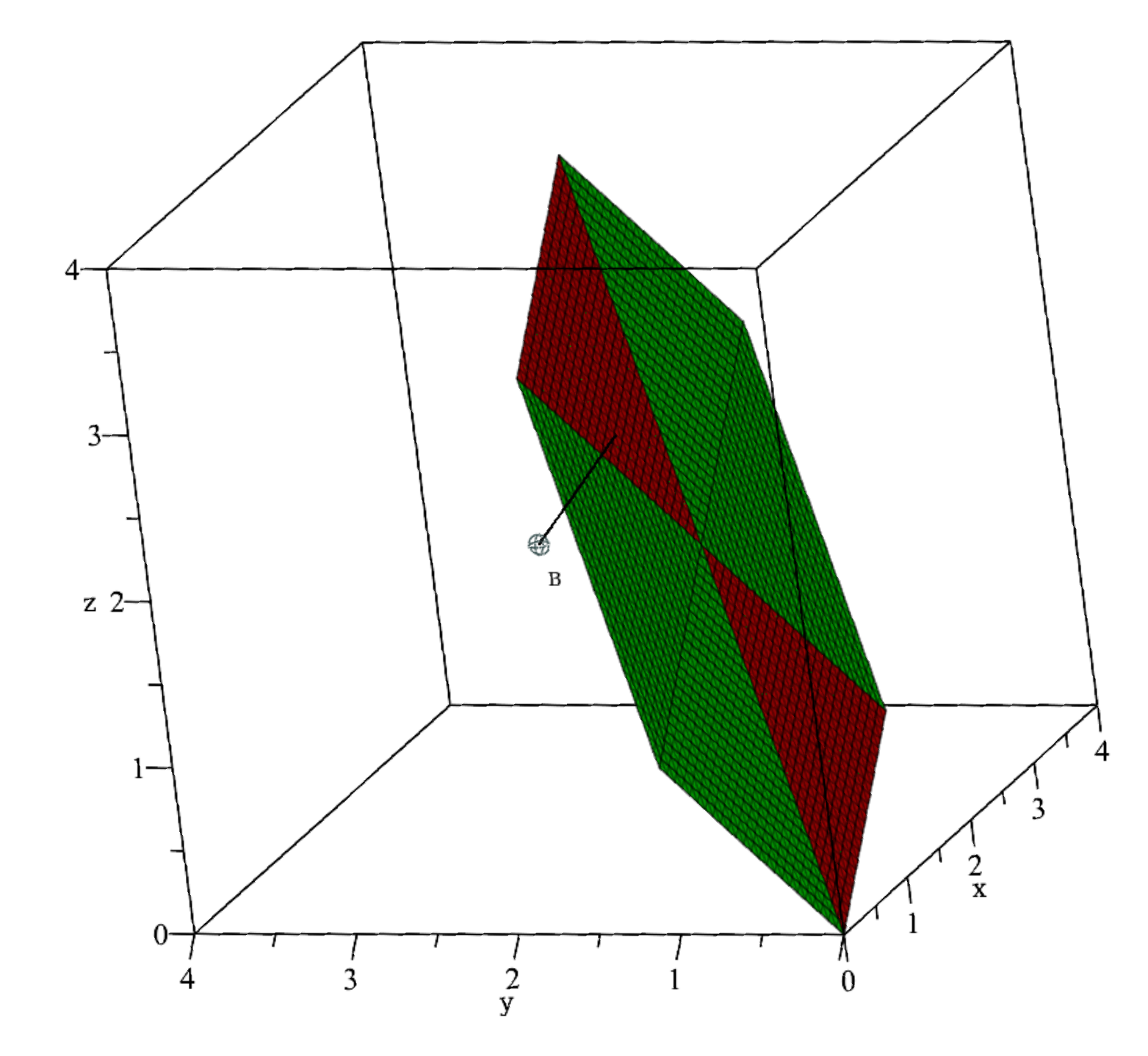}
    \includegraphics[width=3cm]{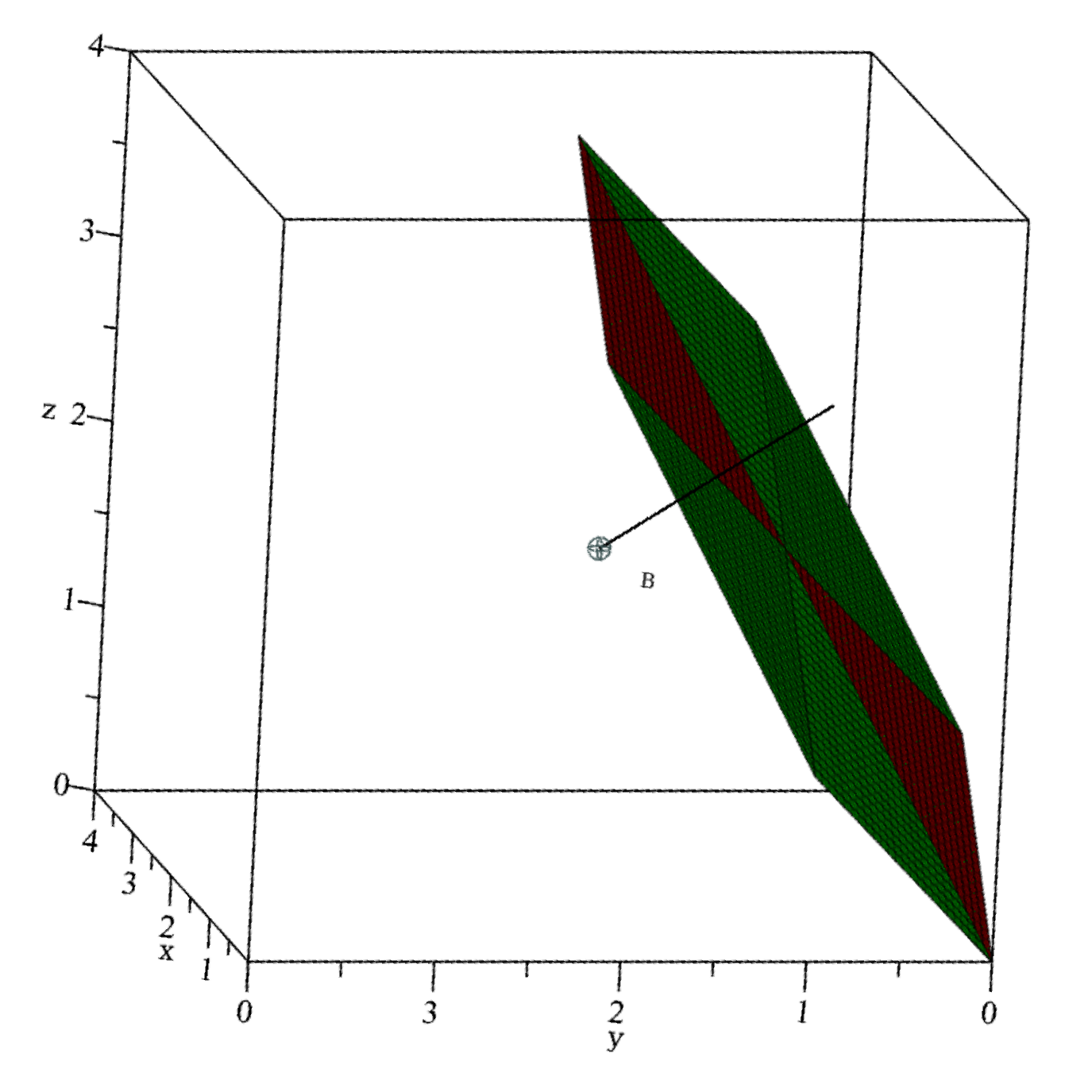}\includegraphics[width=3cm]{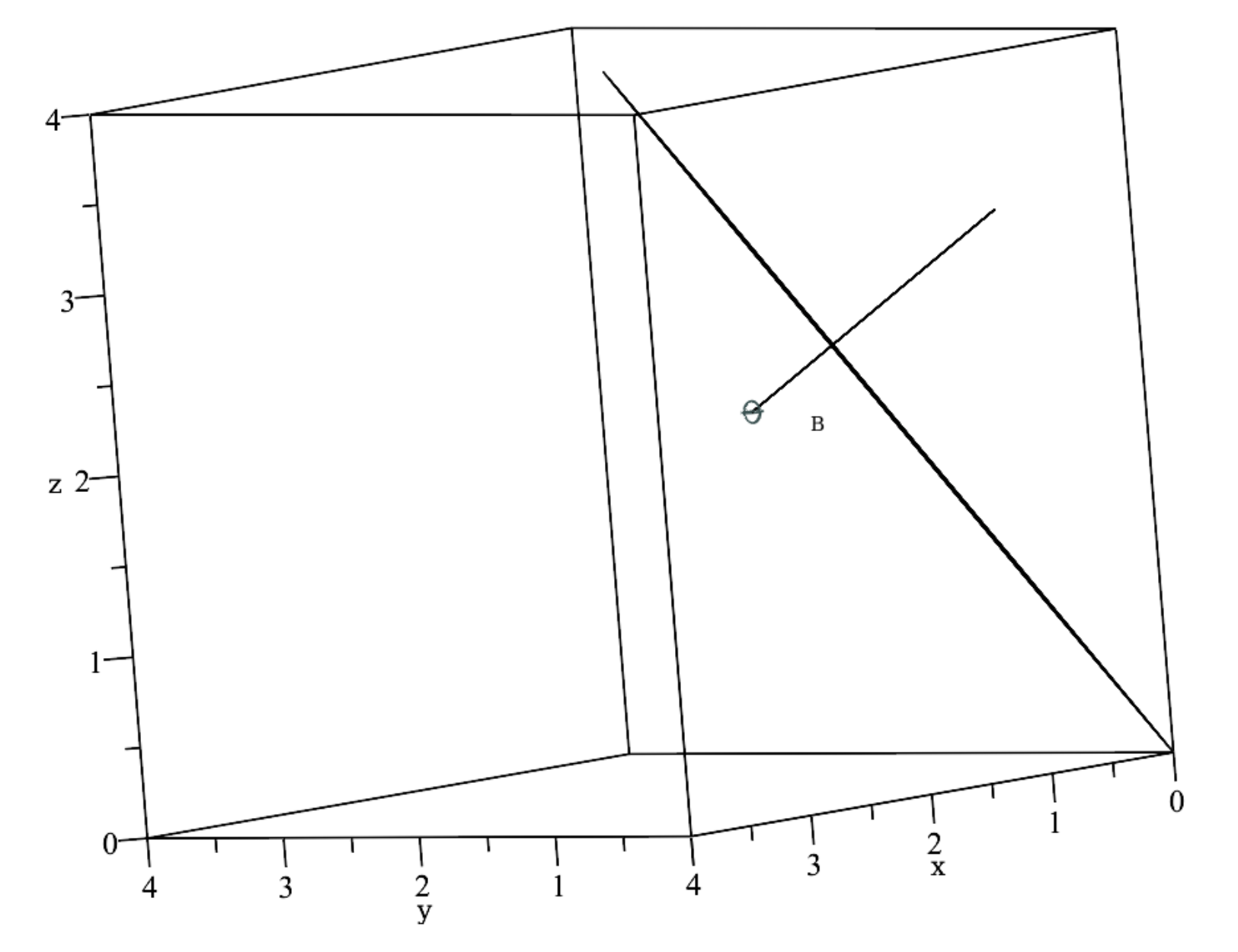}
    \caption{Transformation of the cube with a singular matrix. The three images are an attempt to show in a static medium a student rotating the plot to see that the cube is now flat. $Ax=b$ has no solution because $b$ is not in the plane. We show, however, a projection of $b$ onto the plane, if least-squares is part of the course.  }
    \label{fig:singular}
\end{figure}

An opportunity, neglected in most courses and textbooks, is the making of a connection between equation solving and linear transformations. Typically, a course or book opens with an algebraic account of equation solving. The question of how many solutions an equation has is answered by row reduction and the defining of column space. When transformations are introduced, equation solving is not reconsidered. The equation $Ax=b$ is a transformation of the unknown $x$, in the domain of $A$, to the range, containing $b$. The reverse journey is equation solving, and can be the subject of visualization. In 2-D, everything is rather trivial\footnote{We resisted the temptation to call it ``$2$" trivial.}, so software allowing 3-D interactive plotting is much better. Transforming a cube using
%various non-singular matrices, as shown in Figure \ref{fig:regulartransform}, we explain that every $b$ in the range can be back-tracked to its source $x$ in the domain. The ability of students to  handle, digitally,  the shapes is an important help to their comprehension.
%
% \begin{figure}
%     \centering
%     \includegraphics[width=2cm]{BasicCube.pdf}\includegraphics[width=3cm]{Stretch.pdf}
%     \includegraphics[width=3cm]{Shear.pdf}\includegraphics[width=3cm]{rotate.pdf}
%     \caption{The cube, on the left in the domain of the transformation, is transformed into various shapes, shown as the succession of images to the right: stretch, shear, rotate. }
%     \label{fig:regulartransform}
% \end{figure}
%
%If we transform with
a singular matrix, we observe that the cube is squashed flat. An equation, or the reverse transform, is solvable only if $b$ lies in the plane. See figure \ref{fig:singular}.

\subsection{Programming}
One of the most venerable introductory programming tasks is to write code for LU factoring.  One can then add partial pivoting, complete pivoting, or rook pivoting.  The topic is accessible, but difficult enough that students will really feel a sense of accomplishment when they have succeeded.

The hard part is to get them \textsl{actually to do it} and not to copy someone else's code.  This is especially true in engineering classes, where the students are so heavily pressured that they feel that they \textsl{must} cut corners wherever they can.  One needs to be creative, here, in finding ways to encourage them not to cheat themselves.

One method that we have found effective is to allow them to work in small groups, and to allow them to use code that they find on the internet or copy from other groups \textsl{provided} that they give proper credit and cite where they found it.  Students are frequently surprised that their instructors know about Stack Overflow or Chegg as well; but then, in a work environment, any and all tools will be allowed.  With some creativity in problem assignment, enough novel features can be used so that the online resources will only help, not solve the complete problem for them.  That's unless they use the outright cheating resources where the students post the problems and pay other people to give them the solutions, of course.  To combat that, you have to encourage a culture of honesty by being honest yourself and by actually punishing people caught cheating in that way, so that the honest students feel that they can benefit more by remaining honest.  However, that's a very hard problem to deal with.

It is however something that people in the mathematics mechanization community need to be aware of.  For some decades now, some fully automatic servers have been giving step-by-step solutions to math homework problems.  This is only going to get harder for educators to deal with.  The statement ``if anything can be automated, it should be automated'' ignores the need for the ``White Box'' part of education.  Some concepts need human manual work to be internalized.

\begin{remark}
Many students are only comfortable using computers where they simply enter the data into prescribed fields, and push buttons to achieve pre-programmed aims.  One of the things that we want them to do is to get their ``keyboards dirty'' and engage with a programming language. Doing this at the same time as teaching them the concepts of linear algebra is a stretch.  One should expect only minimal success with getting them to write programs, and then only if you assess them (give them marks) on their ability to do so.  Time spent on that is time that cannot be spent on linear algebra topics.  The topics that we discuss below are chosen in part for their aptness to programming.
\end{remark}
\section{Topics\label{sec:topics}}
In this section we sketch some of the topics that we feel should be encountered in a modern, mechanized, first course in linear algebra, together with how we think that some of the described tools can help with the concepts.
\subsection{The language of matrices}
There is a nontrivial transition from systems of equations such as
\begin{align}
    3x + 4y &= 7 \nonumber\\
    2x - 8y &= 1
\end{align}
to the equivalent matrix equations, and most mechanized systems do not have features to help with this transition.  Matlab, for instance, expects the user to enter the matrices.  We spend some time on this transition, and the conventions that lead to the natural rules for matrix-vector multiplication and thence to matrix-matrix multiplication.  The use of elementary matrices to encode operations on equations (especially elimination of a variable) is a crucial feature.

With beginning students, this takes time.  Hand manipulation is best for this at the beginning, but after experiencing a certain amount of tedium, the students begin to appreciate the ability to construct and manipulate equations through the algebraic rules of matrix multiplication\footnote{They quite like Maple's \lstinline{GenerateMatrix} command, which transforms linear equations with variables into matrix-vector equations.  We try to be careful to introduce this only after the students have some experience in doing the transformation by hand.}.  The simple syntax of Matlab is likely the most appreciated: \lstinline{A*b} for matrix-vector multiplication is close to $\mathbf{A}\cdot \mathbf{b}$, a common human notation; omitting the $\cdot$ seems natural.  Maple's \lstinline{A.b} is somewhat less natural.

Python's notation is similar, except for one thing.
The issue is \textsl{transpose}.  Some linear algebra approaches are very snobbish, and insist that there is no such thing as a row vector or column vector, only abstract vectors.  Python is like this.  This can be very confusing for students.  We have found it best to be explicit and consistent about dimensions in our teaching, and to treat vectors normally as column vectors and to treat these as basically indistinguishable from $n \times 1$ matrices (even that convention needs to be taught: one of our colleagues memorably put it as ``you \textsl{row} with {columns} (oars) when you row a boat'').

The ``four ways'' of interpreting matrix-matrix multiplication is something we explicitly teach.  For instance, in one of these four ways, the matrix-matrix product $\mathbf{A}\mathbf{B}$ can be usefully thought of by first thinking of $\mat{B} = [\mat{b}_1,\mat{b}_2,\ldots,\mat{b}_n]$ as a collection of columns, and then $\mat{A}\mat{B} = [\mat{A}\mat{b}_1,\mat{A}\mat{b}_2,\ldots,\mat{A}\mat{b}_n]$ is then a collection of the column vectors $\mat{A}\mat{b}_k$.  Technological support for this can be as simple as asking the students to construct the matrix on the right hand side explicitly, and verifying that the internal matrix multiplication routine produces the same result.  An advanced question is to consider parallelism in matrix-matrix multiplication using this partition.

We also \textsl{begin} with complex numbers.  They will be needed, so we introduce them first thing.  Without technological support, students hate complex numbers.  With technological support, complex numbers become routine.
\subsection{Parametric Linear Algebra}
One important feature of our course is that it is not purely numerical.  Mathematical modelling frequently involves unknown parameters. One wants the solution in terms of those parameters (if possible) to make it possible to identify those parameters by comparing to experimental data.  There is also the pedagogical value of strengthening student's understanding of formulas, when the answers are not numbers but instead are formulas.

As is well-known in the computer algebra community, this can make computations much more costly and indeed some problems are known to have exponential cost or, worse, combinatorial cost. There is significant literature on the topic, starting with~\cite{sit1992algorithm}.  Recent work includes~\cite{corless2017jordan,CamargosCoutoetal2020a,corless2020parametric} and~\cite{deng2021parametric}.  We will address this issue as it comes up in the various topics. The paper~\cite{deng2021parametric} raises the important point that for many practical problems with only a few parameters, perhaps only one or two, and for problems with structure or low dimension or both, solutions are perfectly feasible using modern computers and infrastructure.

\subsection{Factoring Matrices}
Factoring matrices, whether it is the Turing factoring $\mat{P}\mat{A} = \mat{L}\mat{D}\mat{U}\mat{R}$ which gives the reduced row echelon form~\cite{corless1997turing}, or $\mat{A} = \mat{Q}\mat{R}$ into an orthogonal factor $\mat{Q}$ and upper (right) triangular factor $\mat{R}$, or any of several other factorings, is fundamental for modern linear algebra.  There is the Schur factoring $\mat{A} = \mat{Q} \mat{T}\mat{Q}^H$ which gives the eigenvalues in a numerically stable way.
% \footnote{We remark that most implementations of the Schur factoring today are not satisfactory for teaching, in that this factoring \textsl{should} be continuous in the elements of $\mat{A}$; but is not generally implemented in that way.  This means that a small change in $\mat{A}$ can result in a vastly different pair of matrices $\mat{Q}$ and $\mat{T}$, chiefly re-ordering of eigenvalues $t_{ii}$ or change of signs of entries of $\mat{Q}$.  This can be quite annoying. If someone reading this would be so kind as to implement a continuous version of the Schur factoring, we would be grateful.}.

We teach the notion of factoring matrices as a method of solving linear systems of equations (and of eigenvalue problems).  This represents a conceptual advance over Gaussian Elimination, and has several important consequences in a symbolic context~\cite{jeffrey2013linear,corless1997turing}.  The most important feature in a symbolic context is that a factoring preserves special cases.

Students can factor matrices by hand (and in the beginning, they should).  This gives them something useful to \textsl{do}.  Elementary matrices encoding row operations, column operations, and row exchanges are all useful to teach because they consolidate students' knowledge into a modern framework of understanding of linear algebra, and they do so in a way that allows the student to be \textsl{active}.

Then one can introduce \textsl{block} matrix manipulation and \textsl{block} factoring, with noncommuting elements.  This gives the Schur complement and the Schur determinantal formula.

Interestingly, Maple has recently begun to support matrices over noncommuting variables via the Physics package by Edgardo Cheb--Terrab. This allows students to manipulate block matrices with technology, although they still have to think about dimensions.  This is apparently also possible in SageMath.  Here is an example, showing the Schur complement, in Maple.
\mapleinput
{$ \displaystyle \texttt{>\,} \mathit{with} (\mathit{Physics})\colon  $}

\mapleinput
{$ \displaystyle \texttt{>\,} \mathit{Setup} (\mathit{mathematicalnotation} =\mathit{true})\colon $}

% % \mapleresult
% \begin{dmath}\label{(1)}
% \left[\mathit{mathematicalnotation} =\mathit{true} \right]
% \end{dmath}
\mapleinput
{$ \displaystyle \texttt{>\,} \mathit{Setup} (\mathit{noncommutativeprefix} =\{B \})\colon $}

% \mapleresult
% \begin{dmath}\label{(2)}
% \left[\mathit{noncommutativeprefix} =\left\{\textcolor{olive}{B}\right\}\right]
% \end{dmath}
\mapleinput
{$ \displaystyle \texttt{>\,} \mathit{with} (\mathit{LinearAlgebra})\colon  $}

\mapleinput
{$ \displaystyle \texttt{>\,} A \coloneqq \mathit{Matrix} ([[B [1,1],B [1,2]],[B [2,1],B [2,2]]]\,) $}

% \mapleresult
\begin{dmath}\label{(3)}
A \coloneqq \left[\begin{array}{cc}
\textcolor{olive}{B}_{1,1} & \textcolor{olive}{B}_{1,2}
\\
 \textcolor{olive}{B}_{2,1} & \textcolor{olive}{B}_{2,2}
\end{array}\right]
\end{dmath}
\mapleinput
{$ \displaystyle \texttt{>\,} L \coloneqq \mathit{Matrix} ([[1,0],[B [2,1]\cdot B [1,1]^{-1},1]]) $}

% \mapleresult
\begin{dmath}\label{(4)}
L \coloneqq \left[\begin{array}{cc}
1 & 0
\\
 \textcolor{olive}{B}_{2,1} {\textcolor{olive}{B}_{1,1}}^{-\mathrm{1}} & 1
\end{array}\right]
\end{dmath}
\mapleinput
{$ \displaystyle \texttt{>\,} U \coloneqq \mathit{Matrix} ([[B [1,1],B [1,2]],[0,B [2,2]-B [2,1]\cdot B [1,1]^{-1}\cdot B [1,2]]]\,) $}

% \mapleresult
\begin{dmath}\label{(5)}
U \coloneqq \left[\begin{array}{cc}
\textcolor{olive}{B}_{1,1} & \textcolor{olive}{B}_{1,2}
\\
 0 & \textcolor{olive}{B}_{2,2}-\textcolor{olive}{B}_{2,1} {\textcolor{olive}{B}_{1,1}}^{-\mathrm{1}} \textcolor{olive}{B}_{1,2}
\end{array}\right]
\end{dmath}
\mapleinput
{$ \displaystyle \texttt{>\,} L \cdot U  $}

% \mapleresult
\begin{dmath}\label{(6)}
\left[\begin{array}{cc}
\textcolor{olive}{B}_{1,1} & \textcolor{olive}{B}_{1,2}
\\
 \textcolor{olive}{B}_{2,1} & \textcolor{olive}{B}_{2,2}
\end{array}\right]
\end{dmath}

This \textsl{illustrative} usage of simple noncommuting scalar variables to represent blocks inside matrices, where $1$ represents an appropriately-sized identity matrix and $0$ represents a zero block, might disconcert people intent on formalizing the computations involved. One of the things that would be necessary to properly formalize this would be a notion of dimension of each block; in practice one would want the dimensions to be \textsl{symbolic} but to match appropriately.  We are not aware of any widely-available system at present that can deal properly with this, although there has been research in the area, such as~\cite{Sexton2006,sexton2009computing}.
Making a package widely available that could do such computations correctly would be very welcome.

\subsection{Determinant}
Approaching linear algebra via the determinant is a historically valid approach.  It is pedagogically valid, also, because the students are happier (and better off) with having something to \textsl{do}, not just think about. We feel that it is ``fair game'' that the students be required to memorize the formulas for the determinant and the inverse of a $2\times 2$ matrix (and in fact this memorization is surprisingly useful for them, later).  Laplace expansion (determinant by minors) can be costly and numerically dubious but is extremely useful for sparse symbolic matrices.  More, it is crucial in the one ``gem'' proof that we include in the course simply because it is so pretty, namely the proof of Cramer's Rule\footnote{One of us teaches Cramer's Rule only because of this beautiful proof.  Cramer's Rule itself is not particularly useful computationally nowadays, except in very special situations.  But that proof is so beautiful.  The students seem to like it, too.
%See the appendix to this paper.
} which we learned from~\cite{carlson1992}.

Asking them to memorize a formula for a three-by-three determinant serves no useful purpose, in our opinion, and letting them use technology for computation of third or higher-order determinants seems perfectly justified.

We also demonstrate combinatorial growth by showing the determinant of fully symbolic matrices, for a few small dimensions.  Asking them to program Laplace expansion recursively is also useful for this.  One can also ask them to program the recursive computation of determinant by the Schur determinantal formula $\det\mat{A} = \det \mat{B}_{11} \det( \mat{B}_{22} - \mat{B}_{21}\mat{B}_{11}^{-1}\mat{B}_{12} ) $.  Explicit computation of the inverse of $\mat{B}_{11}$ should be avoided, and can be, by using a suitable factoring.  The end result can be significantly more efficient than Laplace expansion.

We spend time on the geometry of determinant and its relationship to how area transforms under linear transformations; this is needed in calculus, and can be motivating for the students as well because it makes a connection to something that they already know.  Computer visualizations help, here.  The ones freely available on YouTube, especially the very professionally produced ones by \verb+3Blue1Brown+ such as \url{https://youtu.be/Ip3X9LOh2dk}, are hard to compete with.  So, we do not compete, and instead share our favourites (such as that one) with the student.

With determinant in hand, the students have a worthwhile test for linear dependence.  We extend this using the SVD because in the context of data error (which our clientele will surely encounter), the notion of exact singularity or dependence is less useful than that of ill-conditioning or near-dependence.
\paragraph{Least squares}
Matlab will silently return a least-squares solution to overdetermined problems.  Or, even, inconsistent problems.  Therefore it is incumbent on us as instructors to teach least squares solutions, in order that the user may understand and appreciate what the system has done.
\subsection{Eigenvalues and floating-point}
\begin{figure}
    \centering
    \includegraphics[width=0.8\textwidth]{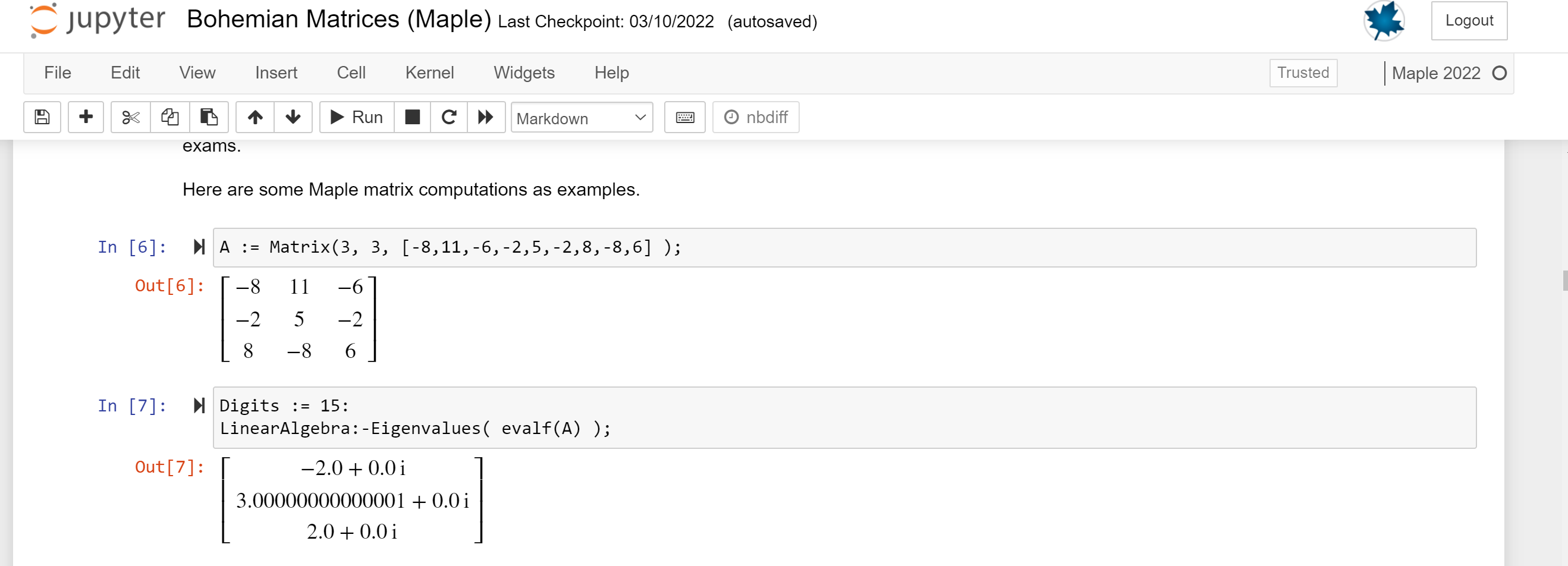}
    \caption{Using Maple from a Jupyter notebook}
    \label{fig:JupyterMapleEigenvalue}
\end{figure}

We teach eigenvalues more by the ``Black Box'' / ``White Box'' approach, because computing eigenvalues by first computing the determinant of $\lambda\mat{I} - \mat{A}$ and then solving the polynomial is a pretty brutal hand computation for anything more than $2 \times 2$ matrices.  We show them what eigenvalues and eigenvectors \textsl{are} by the use of \lstinline{eigshow} or similar, and then set them to compute eigenvalues by the technology.  For instance in Figure~\ref{fig:JupyterMapleEigenvalue} we see how to do this using Maple (from inside a Jupyter notebook).  This requires a discussion of floating-point arithmetic and backward error analysis, which we do not shy away from.  Again, our clientele will encounter data error and they must learn tools such as the condition number (which is really just the derivative) to deal with it; putting numerical error on the same footing as data error gives them the tools to deal with that, as well.  The computation of eigenvalues of small matrices (say, of dimension less than 1000) is a solved problem nowadays.

Indeed we view eigenvalues as answers nowadays because the algorithms are so good in practice (and have recently been shown to be globally convergent in theory, as well~\cite{banks2022global}).  We have had units (in some of our courses) where we talk about companion matrices of various kinds, as tools for solving polynomial equations and systems of polynomial equations.  We discuss this in section~\ref{sec:companions}.

Eigenvalues of \textsl{parametric} matrices are important, for instance in dynamical systems, and their study leads directly to bifurcation theory.  We do not include many such problems, but we have used one in particular, namely a perturbation of Matlab's \lstinline{gallery(3)} matrix to examine the sensitivity of its eigenvalues to perturbations.  This is an advanced topic, however, and occurs only toward the end of the first course (and much more frequently in the second or later course).
\subsection{Special matrices}
There are countless kinds of special matrices.  Likely the most important in practice are symmetric (Hermitian) positive definite matrices; others include orthogonal (unitary) matrices, triangular matrices, banded matrices, circulant matrices, Toeplitz matrices, Hankel matrices, and totally positive matrices.  Getting the students to write programs that generate some of these, or factor some of these in special ways, is quite interesting.  The Cayley transform is quite important nowadays (see e.g.~\cite{mehrmann1996step}) in control theory and in some kinds of scientific computing, and getting students to parameterize orthogonal matrices using symmetric matrices and the Cayley transform may teach several lessons.

While this course should include some of the most common and useful kinds of special matrices, we feel it is also important to let the students invent some of their own kinds of matrices.  Examples of student-generated matrices include ``checkerboard'' matrices which alternate nonzero entries with zero entries and ``anti-tridiagonal'' matrices.  We have found it fun to let the students play, as they program.  Sometimes even their bugs give rise to interesting developments.
\paragraph{Symmetric Positive Definite matrices}
\begin{quote}
    ``Symmetric positive definiteness is one of the highest accolades to which a matrix can aspire.''
\end{quote}
\begin{flushright}
---Nicholas J.~Higham, in~\cite[p.~196]{higham02:ASN}
\end{flushright}
Symmetric Positive Definite (SPD) matrices arise very often in practice.  For an enlightening discussion of just why this is so, see~\cite{Strang1986}.
% In addition to several useful theoretical properties, such as the defining property $\mat{A}$ is SPD if and only if $\mat{x}^T \mat{A} \mat{x} > 0$ for nonzero vectors $\mat{x}$, there are several computational implications, especially that pivoting is not needed and so LU factoring can be carried out efficiently.  The Cholesky factoring $\mat{A} = \mat{L}^T\mat{L}$ is unique, although for symbolic computation the equivalent factoring with \textsl{unit} lower triangular matrices $\mat{L}$ and a positive diagonal matrix $\mat{D}$, namely $\mat{A} = \mat{L}^T\mat{D}\mat{L}$, is frequently better.
The inductive proof of unicity of the Cholesky factoring for SPD matrices (see e.g.~\cite[p.~196]{higham02:ASN}) can be turned into a recursive program for its computation, and this is a useful programming exercise for the students.
The many applications of SPD matrices can be motivating for students, but having the technology to solve them is clearly essential.
\paragraph{Companion matrices\label{sec:companions}}
\begin{quote}
    ``What does this all have to do with matrices? The connection is through the \textsl{companion} matrix.''
\end{quote}
\begin{flushright}
---Cleve Moler, in~\cite{moler1991}.
\end{flushright}
Another thing technology really makes possible is the use of companion matrices and resultants in the solution of polynomial equations.  The topic is surprisingly rich, not just useful.  Algebraically, companion matrices for a monic polynomial $p(z)$ are matrix representations of multiplication by $z$ in the ideal generated by $p(z)$.  Companion matrices are not unique, and indeed there are open problems as to which is the ``best'' companion for a given polynomial $p(z)$, as we will discuss.  Extending the idea to non-monic polynomials leads to generalized eigenvalue problems $p(z) = \det( \lambda \mat{B} - \mat{A} )$ where now $\mat{B}$ is not necessarily the identity matrix (or of full rank).  Using other polynomial bases (e.g. Chebyshev, Bernstein, or Lagrange interpolational bases) leads again to surprisingly deep waters.  Given a (monic) polynomial over the integers, one can ask which companion matrix over the integers has minimal height?  The ``height'' of a matrix is the infinity norm of the matrix made into a vector; that is, the largest absolute value of any entry.  No good algorithms for this problem are known~\cite{chan2019minimal}.  In the case of Mandelbrot polynomials $p_0 = 0$ and $p_{n+1}(z) = z p_n^2(z) +1$ there are companions of height $1$, while the maximum coefficient of $p_n(z)$ is exponential in the degree of $p_n(z)$ (and therefore doubly exponential in $n$).  Smaller height matrices seem to be easier to compute eigenvalues for.

\subsection{Proof and formal methods}
\begin{quote}
``I have absolutely no interest in proving things that I know are true.''
\end{quote}
\begin{flushright}
---the American physicist Henry Abarbanel, at a conference in 1994
\end{flushright}

Entering students in North America have long since been deprived of an introductory course on proof (which was, classically, Euclidean geometry). Typically the first course in which they encounter ``proof'' nowadays is their first linear algebra class.  For the clientele described previously, we feel it is more important to \textsl{motivate} proof at this stage.  Students who are asked to listen to a proof of something they consider obvious (or for which they would be happy to take the professor's word, such as $\det \mat{A}\mat{B} = \det\mat{A} \cdot \det\mat{B}$) do not learn much.  Ed Barbeau put it thus: ``there should be no proof without doubt'' (on the part of the student).

Asking students to write programs is, we believe, a useful intermediate step.  In addition to developing the necessary habit of precise thinking, writing programs makes students receptive to the idea of proving their programs correct (after they have witnessed a few failures, which are somehow always surprising to beginning programmers).

\section{Assessment\label{sec:assessment}}
Assessment is critical for the success of a course.  Students want bribes (marks) in order to spend time on any particular topic.  If a topic is not assessed, then it can be safely skipped and the student can rationally spend their effort on topics that actually will be assessed.

The recent introduction of chat AIs that generate plausible-sounding answers has thrown a further monkey wrench into assessment of courses by project, a method that we have heretofore favoured.  It is even the case that these chat AIs can, perhaps by plagiarising GitHub and other software sources, provide readable (and sometimes even working) software to students.  We may have to go back to individual exams with direct supervision: essentially, oral examinations.  This is so labour intensive that it seems impractical for the very large linear algebra classes that our Universities want us to teach, however.

There are several strategies for written exams that still may be of interest, however, and we give some of them here.  The first is the venerable multiple-choice exam.  For computations, one can remove the ``reverse engineering'' method by asking not for the exact answer, but rather asking for the closest answer not larger than the true answer.  For instance, supposing that the true answer was $\sqrt{2}$, one could list decimal answers (a) $1.2$ (b) $1.3$ (c) $1.5$ and (d) $1.8$.  The desired answer would be (b), $1.3$ .  This tool is surprisingly effective, although many students view it as being ``unfair.''
% Constructing good multiple-choice exams is a skill that can be taught, and there is significant research on it.

A second assessment strategy is to use computer-generated individual questions, where the student is expected to work at their computer (or at a locked-down lab computer) and provide full notes on their work.  These kinds of exams are very stressful for students, however.  They are even more stressful if intrusively-monitoring software is involved (and there may be human rights abuses committed by those pieces of software which the instructor or administration will be responsible for).

For the purpose of discussion, we will assume that no intrusive monitoring software is used, and that measures are taken to alleviate student stress: for instance, one can give out ``practice'' exams ahead of time.

Since we want to include the use of mechanized tools into the assessment, testing in a computational environment is quite natural.  If the students know that they will be tested on their competence in (say) Python, then they will spend some effort to learn it.  Incorporating personalized questions into such exams then becomes both feasible and informative.

% We like the ``ungrading'' or ``mastery grading'' approach where students are allowed to resubmit earlier work until it is perfected.  Again, this does not scale well with class size because it is labour intensive.  But there is no doubt that students learn better when this strategy is involved.
\section{Promoting agreement on syllabus change\label{sec:reaction}}
Some of our colleagues and administrative structures have been very supportive of innovation along these lines.  Others have been, well, reactionary.  Using technology is more labour intensive than is re-using the same old linear algebra textbooks, problem sets, and exam questions.
Using technology also requires continual re-training because the technologies keep changing.
Some people resent being told that they have to change in order to do their jobs well in a changing environment.

We give an example here of a suboptimal linear algebra exam question, taken from last year's multi-section course at Western\footnote{A simple web search for ``Math 1600 Western'' brings the entire exam up, if you wish to see the entire context.}, taught both by progressive and regressive colleagues.  The exam took place without notes, books, calculators, or computers.  Students are allowed by law (in some parts of the world) to have access to their phones, but many universities will attempt to restrict that, too.  The exams at Western typically have quite alarming language on the cover sheet saying that students caught with a cell phone will be given a zero.  We feel that this is a lamentable state.

The question was: Find the inverse of the matrix
\begin{equation}
    \mat{A} = \left[\begin{array}{rrr}
\phantom{-}2 & \phantom{-}1 & 0
\\
 1 & 0 & -1
\\
 0 & 1 & 1
\end{array}\right]\>.
\end{equation}
% The desired answer is
% \begin{equation}
%     \mat{A}^{-1} = \left[\begin{array}{rrr}
% 1 & -1 & -1
% \\
%  -1 & 2 & 2
% \\
%  1 & -2 & -1
% \end{array}\right]\>.
% \end{equation}
This question does have a few virtues.  For one, it is something the students can \textsl{do}.  It was worth three marks, which the students could grind out.

But it also has some serious flaws.  Probably the most serious is that it does not test anything that the students will really need in their future use of linear algebra.  There were calculators thirty years ago that could solve this problem in under a second.  No one is going to invert $3 \times 3$ matrices by hand any more, unless there is something special about it.  [There \textsl{is} something special about this matrix; it is unimodular, so that the elements of the inverse are all integers.  That didn't happen by chance, so we suspect the examiners chose the question so as not to strain the student's arithmetic overmuch.]

More, not only will students not need to invert by hand, they usually will not need to invert at all.  The \textsl{inversion of matrices} is really only of very specialized concern nowadays.  There are statistical applications where the elements of the inverse are what is wanted; but for the most part, ``Anything that you can do with the matrix inverse can be done without it.'' Matrix factorings are much more important.

Students are rational creatures.  If this is the kind of question that they have to answer in order to pass, then they will spend their time trying to find strategies to give good answers to this kind of question.  They will do that at the expense of time spent learning to program (for instance).

This represents a significant lost opportunity for the student and for this University.
Indeed, the absolute explosion in on-line courses (for instance, at \href{https://brilliant.org}{brilliant.org}, where they claim that interactive learning is six times more effective than lectures) is a direct response to the failure of many universities to adapt their courses.  Students resent having to pay twice to get the knowledge they actually want and need.
The next few years are going to be ``interesting.''

One way to repair that particular question might be to ask if the matrix factors into a lower triangular and upper triangular factor, without pivoting.  The matrix is tridiagonal, so this variation has fewer computations, although this time involving fractions (just $1/2$ though).  This is something that could be asked even if the student has access to technology during the exam.  The details of the computation are not that important---it is just arithmetic---but the question of whether or not the factoring can be done without pivoting would require some understanding of the process involved.
% For interest, here is the answer: yes, it can be done, as follows.
% \begin{equation}
%     \mat{A} = \left[\begin{array}{ccc}
% 1 & 0 & 0
% \\
%  \frac{1}{2} & 1 & 0
% \\
%  0 & -2 & 1
% \end{array}\right]
% \left[\begin{array}{ccc}
% 2 & 1 & 0
% \\
%  0 & -\frac{1}{2} & -1
% \\
%  0 & 0 & -1
% \end{array}\right]\>.
% \end{equation}
% Notice that the factors are bidiagonal.  This factoring exemplifies a more efficient method of solution than using the (full) matrix inverse.
\section{Concluding Remarks}
The state of the art for learning linear algebra is, to our minds, unsatisfactory, though getting better.  Technological platforms are split: some are proprietary, while some others are unsupported at the level needed for reliable use.  Methods and syntax are not standardized (or, rather, there are too many standards).  The textbooks largely do not integrate mechanized mathematical tools into the learning process. [A very notable exception is~\cite{VanLoan2010}, which uses Matlab extensively.]  Yet failing to use a mechanized approach does a true disservice to students who will go on to practice linear algebra in some kind of mechanized environment.

The role of technology, including formal methods, is therefore multiplex.  We believe that people must be trained in its use.  In particular, people must be trained to want proof, and to want formal methods.  We feel that having students write their own programs plays a motivating role in that training as well as a developmental role.  The first linear algebra course is important not only because its tools and concepts are critical for science, but also as a venue for teaching the responsible use of mathematical technology.
\bibliographystyle{splncs04}

\end{document}